\documentclass[11pt, twoside]{article}
\usepackage{latexsym}
\usepackage{amsmath}
\usepackage{amssymb}
\usepackage[all]{xy}
\usepackage{amsfonts}
\usepackage{verbatim}
\usepackage{amsthm}
\usepackage{mathrsfs}
\usepackage{epsfig}
\usepackage{xy}
\usepackage{array}
\usepackage{stmaryrd}
\usepackage{graphicx,color}
\usepackage{xcolor}
\usepackage{tikz}
\usetikzlibrary{arrows,calc}
\usepackage{etex}
\usepackage{mathdots}
\usepackage{float}
\usepackage{graphics}
\usepackage{pdflscape}
\usepackage{mathrsfs}
\newcommand{\add}{\mathsf{add}\hspace{.01in}}
\usepackage{anysize,hyperref}
\input xypic
\xyoption{all}

\usepackage[perpage,symbol]{footmisc}
\usepackage{setspace}
\def\C{\mathscr{C}}
\def\E{\mathbb{E}}
\def\s{\mathfrak{s}}

\def\op{^\mathrm{op}}

\def\del{\delta}
\def\dr{\ar@{->}[r]}

\def\X{\mathscr{X}}\def\S{\mathscr{S}}\def\T{\mathscr{T}}
\def\Y{\mathscr{Y}}
\def\add{\mbox{add}}

\def\Hom{\mbox{Hom}}

\begin{document}
\baselineskip=15pt
\title{\Large{\bf Semi-abelian categories arising from pseudo cluster tilting\\[2mm] subcategories }}
\medskip
\author{Jian He and Jing He\footnote{Corresponding author. ~Jian He is supported by the National Natural Science Foundation of China (Grant No. 12171230) and Youth Science and Technology Foundation of Gansu Provincial (Grant No. 23JRRA825). Jing He is supported by the Hunan Provincial Natural Science Foundation of China (Grant No. 2023JJ40217).}}

\date{}

\maketitle
\def\blue{\color{blue}}
\def\red{\color{red}}

\newtheorem{theorem}{Theorem}[section]
\newtheorem{lemma}[theorem]{Lemma}
\newtheorem{corollary}[theorem]{Corollary}
\newtheorem{proposition}[theorem]{Proposition}
\newtheorem{conjecture}{Conjecture}
\theoremstyle{definition}
\newtheorem{definition}[theorem]{Definition}
\newtheorem{question}[theorem]{Question}
\newtheorem{remark}[theorem]{Remark}
\newtheorem{remark*}[]{Remark}
\newtheorem{example}[theorem]{Example}
\newtheorem{example*}[]{Example}
\newtheorem{condition}[theorem]{Condition}
\newtheorem{condition*}[]{Condition}
\newtheorem{construction}[theorem]{Construction}
\newtheorem{construction*}[]{Construction}

\newtheorem{assumption}[theorem]{Assumption}
\newtheorem{assumption*}[]{Assumption}

\baselineskip=17pt
\parindent=0.5cm

\begin{abstract}
\baselineskip=16pt
The notion of a pseudo cluster tilting subcategory $\X$ in an extriangulated category $\C$ is defined in this article. We prove that the quotient category $\C/\X$, obtained by factoring an extriangulated category by a pseudo cluster tilting subcategory, is a semi-abelian category. Furthermore, we also show that the quotient category $\C/\X$ is an abelian category if and only if certain self-orthogonal conditions are satisfied. As an application, these results generalize the work of Xu and Zheng in the exact category.\\[0.5cm]
\medskip
\end{abstract}

\pagestyle{myheadings}
\markboth{\rightline {\scriptsize   J. He and J. He}}
         {\leftline{\scriptsize   Semi-abelian categories arising from pseudo cluster tilting subcategories}}

\section{Introduction}
Nakaoka and Palu introduced the concept of extriangulated categories in their seminal work \cite{NP}. This notion simultaneously generalizes exact categories and triangulated categories. Exact categories, which include abelian categories, and extension-closed subcategories within extriangulated categories are considered as specific instances of extriangulated categories.
Moreover, there exist additional instances of extriangulated categories that do not fall within the categories of exact or triangulated, as documented in \cite{NP,ZZ1,HZZ}.

Cluster tilting theory offers a method for constructing abelian categories from some triangulated categories and exact categories. Koenig and Zhu \cite{KZ} established a comprehensive framework for transitioning from triangulated categories to abelian categories by factoring out cluster tilting subcategories.
Following this line of research, Demonet and Liu \cite{DL} subsequently demonstrated a similar result to that of Koenig and Zhu. Specifically, they presented a generalized framework for transitioning from exact categories to abelian categories by isolating cluster tilting subcategories.
Expanding on this work, Liu and Nakaoka \cite{LN} demonstrated that any quotient of an extriangulated category modulo a cluster tilting subcategory inherently carries an induced abelian structure. This result extends and generalizes both \cite[Theorem 3.3]{KZ} and
\cite[Theorem 3.2]{DL}.

Let $(\C, \Omega)$ be an exact category with the exact structure $\Omega$. It is well known that the category $\Omega(\C)$ of conflations in $\C$ forms an additive category. Furthermore, the category $\Omega(\C)$ possesses an exact structure, computed degree-wise. Let $S(\C)$ denote the full subcategory of $\Omega(\C)$ consisting of all splitting conflations.
Recently, in their work \cite{XZ}, Xu and Zheng introduced the concept of a pseudo cluster tilting subcategory $\X$ in an exact category $\C$. They demonstrated that $S(\C)$ always constitutes a pseudo cluster tilting subcategory of $\Omega(\C).$ More generally, Xu and Zheng established that the quotient category $\C/\X$ of an exact category by a pseudo cluster tilting subcategory is a semi-abelian category, with $\C/\X$ becoming an abelian category if and only if certain self-orthogonal conditions are satisfied.
Motivated by these findings, a natural question arises: Can the results of Xu and Zheng \cite{XZ} be extended within the framework of extriangulated categories? In this article, we provide an affirmative answer.

Consider an extriangulated category $\C$ with a pseudo cluster tilting subcategory $\X$.
Our first main result establishes that the quotient category $\C/\X$, obtained by dividing an extriangulated category by a pseudo cluster tilting subcategory, becomes a semi-abelian category, as detailed in Theorem \ref{semi-abel}. Our second main result demonstrates that $\C/\X$ attains the status of an abelian category if and only if specific self-orthogonal conditions are satisfied, as outlined in Theorem \ref{abel}. These findings represent a broadening and generalization of the results initially presented by Xu and Zheng in \cite{XZ}.
To further elucidate our main results, we provide illustrative examples.

The paper is organized as follows: In Section 2, we provide a review of essential definitions and facts necessary for our study. In Section 3, we present the proofs of the main results of this article.

\section{Preliminaries}

We briefly recall some definitions and basic properties of extriangulated categories from \cite{NP}.
We omit some details here, but the reader can find them in \cite{NP}.

Let $\C$ be an additive category equipped with an additive bifunctor
$$\mathbb{E}: {\C}^{\rm op}\times \C\rightarrow {\rm Ab},$$
where ${\rm Ab}$ is the category of abelian groups. For any objects $A, C\in\C$, an element $\delta\in \mathbb{E}(C,A)$ is called an $\mathbb{E}$-extension.
Let $\mathfrak{s}$ be a correspondence which associates an equivalence class $$\mathfrak{s}(\delta)=\xymatrix@C=0.8cm{[A\ar[r]^x
 &B\ar[r]^y&C]}$$ to any $\mathbb{E}$-extension $\delta\in\mathbb{E}(C, A)$. This $\mathfrak{s}$ is referred to as a \emph{realization} of $\mathbb{E}$ if it makes the diagrams in \cite[Definition 2.9]{NP} commutative.

A triplet $(\C, \mathbb{E}, \mathfrak{s})$ is called an {\it extriangulated category} if it satisfies the following conditions.
\begin{itemize}
\item $\mathbb{E}\colon\mathcal{C}^{\rm op}\times \mathcal{C}\rightarrow \rm{Ab}$ is an additive bifunctor.

\item $\mathfrak{s}$ is an additive realization of $\mathbb{E}$.

\item $\mathbb{E}$ and $\mathfrak{s}$  satisfy the compatibility conditions $(\rm ET3),(\rm ET3)^{\rm op},(\rm ET4),(\rm ET4)^{\rm op}$ in \cite[Definition 2.12]{NP}.
\end{itemize}

\begin{remark}\label{exam}
We know that both exact categories and triangulated categories are extriangulated categories
\cite[Example 2.13]{NP} and extension-closed subcategories of extriangulated categories are
again extriangulated \cite[Remark 2.18]{NP}. Moreover, there exist extriangulated categories which
are neither exact categories nor triangulated categories \cite[Proposition 3.30]{NP} and \cite[Example 4.14]{ZZ1}.
\end{remark}

We collect the following terminology from \cite{NP}.

\begin{definition}
Let $(\C,\E,\s)$ be an extriangulated category.
\begin{itemize}
\item[(1)] If a conflation $A\overset{x}{\longrightarrow}B\overset{y}{\longrightarrow}C$ realizes $\delta\in\mathbb{E}(C,A)$, we call the pair $(A\overset{x}{\longrightarrow}B\overset{y}{\longrightarrow}C,\delta)$ an {\it $\E$-triangle}, and write it in the following way.
$$A\overset{x}{\longrightarrow}B\overset{y}{\longrightarrow}C\overset{\delta}{\dashrightarrow}$$

\item[(2)] Let $A\overset{x}{\longrightarrow}B\overset{y}{\longrightarrow}C\overset{\delta}{\dashrightarrow}$ and $A^{\prime}\overset{x^{\prime}}{\longrightarrow}B^{\prime}\overset{y^{\prime}}{\longrightarrow}C^{\prime}\overset{\delta^{\prime}}{\dashrightarrow}$ be any pair of $\E$-triangles. If a triplet $(a,b,c)$ realizes $(a,c)\colon\delta\to\delta^{\prime}$, then we write it as
$$\xymatrix{
A \ar[r]^x \ar[d]^a & B\ar[r]^y \ar[d]^{b} & C\ar@{-->}[r]^{\del}\ar[d]^c&\\
A'\ar[r]^{x'} & B' \ar[r]^{y'} & C'\ar@{-->}[r]^{\del'} &}$$
and call $(a,b,c)$ a {\it morphism of $\E$-triangles}.

\end{itemize}
\end{definition}

The following lemma is presented in \cite[Proposition 1.20]{LN}, which is an alternative version of \cite[Corollary 3.16]{NP}.

\begin{lemma}\label{1.20}
 Let $A\overset{x}{\longrightarrow}B\overset{y}{\longrightarrow}C\overset{\delta}{\dashrightarrow}$ be a $\E$-triangle. let $f\colon A \to D$ be any morphism, and let  $D\overset{d}{\longrightarrow}E\overset{c}{\longrightarrow}C\overset{f^{\ast}\delta}{\dashrightarrow}$ be any $\E$-triangle realizing $f^{\ast}\delta$. Then there is a morphism $g$ which gives a morphism of $\E$-triangles
$$\xymatrix{
A \ar[r]^x \ar[d]^f & B\ar[r]^y \ar[d]^{g} & C\ar@{-->}[r]^{\del}\ar@{=}[d]&\\
D\ar[r]^{d} & E \ar[r]^{e} & C\ar@{-->}[r]^{f^{\ast}\del} &}$$
and moreover, $A\overset{\binom{-f}{x}}{\longrightarrow}D\oplus B\overset{(d,g)}{\longrightarrow}E\overset{e^{\ast}\delta}{\dashrightarrow}$ be an $\E$-triangle.

\begin{lemma}\label{exact} \rm{\cite[Proposition 3.3]{NP} }
Let $\xymatrix{A\ar[r]&B\ar[r]&C\ar@{-->}[r]&}$
be an $\E$-triangle in $\C$. Then we have the following two exact sequences:
$$\C(-, A)\xrightarrow{~}\C(-, B)\xrightarrow{~}\C(-, C)\xrightarrow{~}
\E(-, A)\xrightarrow{~}\E(-, B);
$$
$$\C(C,-)\xrightarrow{~}\C(B,-)\xrightarrow{~}\C(A,-)\xrightarrow{~}
\E(C,-)\xrightarrow{~}\E(B,-).$$
\end{lemma}

\end{lemma}
\section{Quotient categories of extriangulated categories by pseudo cluster tilting subcategories}
In this section, when we say that $\mathcal{T}$ is a subcategory of $\mathcal{C}$, we always mean that $\mathcal{T}$ is a full subcategory.

Let $\C$ be an additive category and $\X$ a subcategory of $\C$.
Recall that we say a morphism $f\colon A \to B$ in $\C$ is an $\X$-\emph{monic} if
$$\C(f,X)\colon \C(B,X) \to \C(A,X)$$
is an epimorphism for all $X\in\X$. We say that $f$ is an $\X$-\emph{epic} if
$$\C(X,f)\colon \C(X,A) \to \C(X,B)$$
is an epimorphism for all $X\in\X$.
Similarly,
we say that $f$ is a left $\X$-approximation of $A$ if $f$ is an $\X$-monic and $B\in\X$.
We say that $f$ is a right $\X$-approximation of $B$ if $f$ is an $\X$-epic and $A\in\X$.

We denote by $\C/\X$
the category whose objects are objects of $\C$ and whose morphisms are elements of
$\Hom_{\C}(A,B)/\X(A,B)$ for $A,B\in\C$, where $\X(A,B)$ is the subgroup of $\Hom_{\C}(A,B)$ consisting of morphisms
which factor through an object in $\X$.
Such category is called the (additive) quotient category
of $\C$ by $\X$. For any object $C\in \C$, we denote by $\underline{C}$ the image of $C$ under
the natural quotient functor $\C\to\C/\X$. For any morphism $f\colon A\to B$ in $\C$, we denote by $\underline{f}$ the image of $f$ under
the natural quotient functor $\C\to\C/\X$.

\begin{definition}\label{dd1}\rm{\cite[Definition 3.21]{ZZ1}}
Let $\mathcal{C}$ be an  extriangulated category. A subcategory $\mathcal{T}$ of $\mathcal{C}$ is called
\emph{strongly contravariantly finite}, if for any object $C\in\mathcal{C}$, there exists an $\E$-triangle
$$\xymatrix{K\ar[r]&T\ar[r]^{g}&C\ar@{-->}[r]^{\del}&,}$$
where $g$ is a right $\mathcal{T}$-approximation of $C$.

Dually, a subcategory $\mathcal{T}$ of $\mathcal{C}$ is called
\emph{strongly  covariantly  finite}, if for any object $C\in\mathcal{C}$, there exists an $\E$-triangle
$$\xymatrix{C\ar[r]^{f}&T\ar[r]&L\ar@{-->}[r]^{\del'}&,}$$
where $f$ is a left $\mathcal{T}$-approximation of $C$.

A strongly contravariantly finite and strongly  covariantly finite subcategory is called \emph{ strongly functorially finite}.
\end{definition}

\begin{definition}\cite[Definition 2.10]{ZZ1}
Let $(\C, \E, \s)$ be an extriangulated category and $\X$ a subcategory of $\C$.
\begin{itemize}
\item $\X$ is called \emph{rigid} if $\E(\X,\X)=0$;
\item $\X$ is called \emph{cluster tilting} if it satisfies the following conditions:
\begin{enumerate}
\item[(1)] $\X$ is a strongly functorially finite in $\C$;
\item[(2)] $M\in \X$ if and only if $\E(M,\X)=0$;
\item[(3)] $M\in \X$ if and only if $\E(\X,M)=0$.
\end{enumerate}
\end{itemize}
\end{definition}

By definition of a cluster tilting subcategory, we can immediately conclude:
\begin{remark}\cite[Definition 2.10]{LZ}\label{y6}
Let $(\C, \E, \s)$ be an extriangulated category  and $\X$ a subcategory of $\C$.
\begin{itemize}
\item  $\X$ is a cluster tilting subcategory of $\C$ if and only if
\begin{enumerate}
\item[$(1^{'})$] $\X$ is rigid;
\item[$(2^{'})$] For any $C\in\C$, there exists an $\E$-triangle $\xymatrix@C=0.5cm{C\ar[r]^{a\;} & X_1 \ar[r]^{b} & X_2\ar@{-->}[r]^{\del}&,}$ where $X_1,X_2\in\X$ and $a$ is a left $\X$-approximation of $C$;
\item[$(3^{'})$] For any $C\in\C$, there exists an $\E$-triangle $\xymatrix@C=0.5cm{X_3\ar[r]^{c\;} & X_4 \ar[r]^{d} & C\ar@{-->}[r]^{\eta}&,}$ where $X_3,X_4\in\X$ and $d$ is a right $\X$-approximation of $C$.
\end{enumerate}
\end{itemize}
\end{remark}

\begin{definition}\label{pct}
Let $(\C, \E, \s)$ be an extriangulated category  and $\X$ a subcategory of $\C$. $\X$ is called a pseudo cluster tilting subcategory of $\C$, if it satisfies the conditions $(2^{'})$ and $(3^{'})$ of Remark \ref{y6}.
\end{definition}

\begin{remark}\label{epst}
If the category $\C$ is exact, then the Definition \ref{pct} coincides with the definition of pseudo cluster tilting subcategory (cf. \cite{XZ}).
\end{remark}

Now we give some examples of pseudo cluster tilting subcategories.

\begin{example}\label{00ex}
 Let $A=kQ/I$ be a self-injective algebra given by the quiver $$Q:
\setlength{\unitlength}{0.03in}\xymatrix{1 \ar@<0.5ex>[r]^{\alpha}_{\ } & 2\ar@<0.5ex>[l]^{\beta}}$$ and $I=\langle\alpha\beta\alpha\beta, \beta\alpha\beta\alpha\rangle$. Let $\C$ be the stable module category
${\rm \underline{\rm mod}}A$ of $A$. This is a triangulated category, and
we describe the Auslander-Reiten quiver of $\C={\rm \underline{\rm mod}}A$ in the following picture:
$$
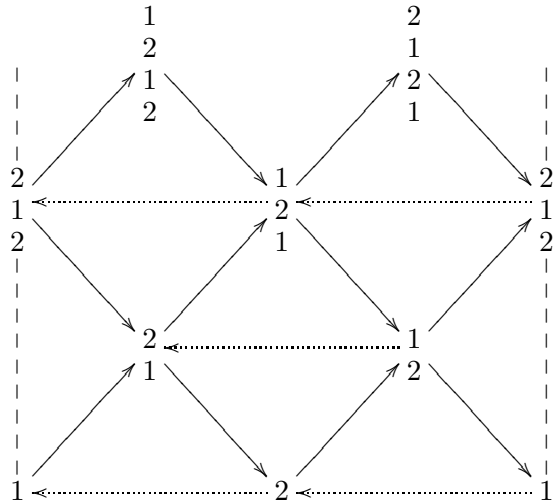
\xymatrix@!@C=0.01cm@R=0.2cm{
&\txt{1\\2\\1\\2}\ar[dr]&&\txt{2\\1\\2\\1}\ar[dr]&\\
{\txt{2\\1\\2}}\ar[ur]\ar[dr]\ar@{--}[u]&&\txt{1\\2\\1}\ar@{.>}[ll]
\ar[ur]\ar[dr]&&\txt{2\\1\\2}\ar@{.>}[ll]\ar@{--}[u]\\
&\txt{2\\1}\ar[dr]\ar[ur]&&\txt{1\\2}\ar[dr]\ar[ur]\ar@{.>}[ll]&\\
 1\ar[ur]\ar@{--}[uu]&&{2}\ar@{.>}[ll]\ar[ur]&&1\ar@{.>}[ll]\ar@{--}[uu] }
$$
where the leftmost and rightmost columns are identified. We denote by $\add(X)$
the subcategory of $\C$ consisting of direct summands of direct sum of finitely many copies
of $X$.

It is straightforward to verify that the subcategory
$$\add(2\oplus\begin{aligned}1\\[-3mm]2\end{aligned}\oplus \begin{aligned}2\\[-3mm]1\\[-3mm]2\end{aligned})$$
is a pseudo cluster tilting subcategory of $\C$, but not a cluster tilting subcategory of $\C$, since ${\rm Hom}_{\C}(2,\begin{aligned}1\\[-3mm]2\end{aligned}[1])=
{\rm Hom}_{\C}(2,\begin{aligned}1\\[-3mm]2\end{aligned})\neq 0$.
\end{example}

\begin{example}\label{ex2}
Let $Q$ be the quiver $\setlength{\unitlength}{0.03in}\xymatrix{1 \ar[r]^{\alpha} & 2}.$  Assume that $\tau_Q$ is the Auslander-Reiten translation in $D^b(kQ)$. We consider the repetitive cluster category $\C=D^b(kQ)/\langle \tau_Q^{-2}[2]\rangle$ introduced by Zhu in \cite{Z}, whose objects are the same in $D^b(kQ)$, and whose morphisms are given by $$\Hom_{D^b(kQ)/\langle \tau_Q^{-2}[2]\rangle}(X, Y)=\bigoplus_{i\in \mathbb{Z}}\Hom_{D^b(kQ)}(X, (\tau_Q^{-2}[2])^iY).$$  It is shown in \cite{Z} that $\C$ is a triangulated category.
We describe the Auslander-Reiten quiver of $\C$ in the following picture:
$$
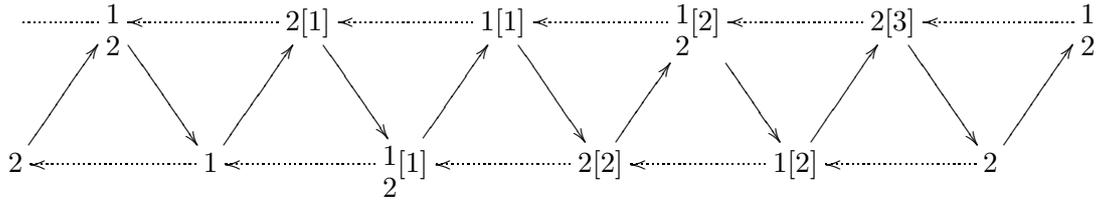
\xymatrix@!@C=0.4cm@R=1cm{  
&\txt{1\\2}\ar@{.}[l]\ar[dr]&&{2[1]}\ar@{.>}[ll]\ar[dr]&&1[1]
\ar@{.>}[ll]\ar[dr]&&{\txt{1\\2}[2]}\ar@{.>}[ll]\ar[dr]&&2[3]\ar@{.>}[ll]\ar[dr]&&\txt{1\\2}\ar@{.>}[ll]  \\
2\ar[ur] &&{1}\ar@{.>}[ll] \ar[ur]&&\txt{1\\2}[1]\ar@{.>}[ll]\ar[ur] && 2[2]\ar@{.>}[ll] \ar[ur]&&{1[2]}\ar@{.>}[ll] \ar[ur]&&2 \ar@{.>}[ll]\ar[ur]  }$$
It is straightforward to verify that the subcategory
$$\add(2\oplus1\oplus\begin{aligned}2\\[-3mm]1\end{aligned}\oplus 2[2]\oplus\begin{aligned}1\\[-3mm]2\end{aligned}[2])$$
is a pseudo cluster tilting subcategory of $\C$, but not a cluster tilting subcategory of $\C$, since ${\rm Hom}_{\C}(1,2[1])\neq 0$.
\end{example}

\begin{example}\label{expst}
This example comes from the reference \cite{XZ}.
As a special  extriangulated category, let $(\C, \Omega)$ be an exact category with the exact structure $\Omega$. We consider the category of conflations of $\C$, denoted by $\Omega(\C)$. It is well known $\Omega(\C)$ is an additive category. The exact structure of $\Omega(\C)$ is the usual exact structure computed degree-wise, written as $(\Omega(\C),\Omega)$. Let $S(\C)$ be the full subcategory of $\Omega(\C)$ consisting of all splitting conflations. Then  $S(\C)$ is a pseudo cluster tilting subcategory of $(\Omega(\C),\Omega)$.
\end{example}

We recall from \cite{B,W} an additive $\C$ with kernels and cokernels is called \emph{semi-abelian} if for any morphism $f\colon X\to Y$, the canonical morphism $\hat{f}\colon {\rm Coim}f\to {\rm Im}f$ is both an epimorphism and a monomorphism. It is abelian if $\hat{f}$ is an isomorphism. Sometimes, A morphism which is both an epimorphism and a monomorphism is called regular. It can be seen that a semi-abelian category is abelian if and only if regular morphisms and isomorphisms coincide.

Our first main result is the following.
\begin{theorem}\label{semi-abel}
Let $\C$ be extriangulated category and $\X$ a pseudo cluster tilting subcategory of $\C$. Then the quotient category $\C/\X$ is a semi-abelian category.
\end{theorem}
In order to prove Theorem \ref{semi-abel}, we need some preparations as follows.

An $\E$-triangle $A\overset{}{\longrightarrow}B\overset{}{\longrightarrow}C\overset{}{\dashrightarrow}$ in $\C$ is called to be $\C(\X,-)$-exact (resp. $\C(-,\X)$-exact), if for any $X\in \X$, the induced sequence of abelian group $\C(X,A)\xrightarrow{}\C(X,B)\xrightarrow{}\C(X,C)\xrightarrow{}0$ (resp. $\C(C,X)\xrightarrow{}\C(B,X)\xrightarrow{}\C(A,X)\xrightarrow{}0$) is exact in ${\rm Ab}$.

An $\E$-triangle $A\overset{x}{\longrightarrow}B\overset{y}{\longrightarrow}C\overset{\delta}{\dashrightarrow}$ is called split if $\delta=0$. By \cite[Corollary 3.5]{NP}, we know that it is split if and only if $x$ is a section or $y$ is a retraction.

\begin{lemma}\label{l1}
Let $(\C, \E, \s)$ be an extriangulated category  and $\X$ a subcategory of $\C$.
\begin{itemize}
\rm\item[(1)] If $\X$ is strongly covariantly finite, then for any morphism $f\colon X\to Y$ in $\C$, there exist some $M\in\X$ and an $\E$-triangle $X\xrightarrow{}Y\oplus M\xrightarrow{}Z\overset{}{\dashrightarrow}$ which is $\C(-,\X)$-exact.
\item[(2)] If $\X$ is strongly contravariantly finite, then for any morphism $g\colon Y\to Z$ in $\C$, there exist some $N\in\X$ and an $\E$-triangle $X\xrightarrow{}Y\oplus N\xrightarrow{~b~}Z\overset{}{\dashrightarrow}$ which is $\C(\X,-)$-exact.

\end{itemize}
\end{lemma}

\proof (1) Since $\X$ is strongly covariantly finite, then for $X\in\C$, there exists an $\E$-triangle $$X\xrightarrow{h}M\xrightarrow{h'}Q\overset{\delta}{\dashrightarrow},$$ where $h$ is a left $\X$-approximation of $X$. By Lemma \ref{1.20}, we have the following commutative diagram
$$\xymatrix{
X \ar[r]^h \ar[d]^f & M\ar[r]^{h'} \ar[d]^{l} & Q\ar@{-->}[r]^{\del}\ar@{=}[d]&\\
Y\ar[r]^{k'} & Z \ar[r]^{k'} & Q\ar@{-->}[r]^{f_{\ast}\del} &}$$
of $\E$-triangles in $\C$, and moreover, $X\overset{\binom{-f}{h}}{\longrightarrow}Y\oplus M\overset{(k',l)}{\longrightarrow}Z\overset{}{\dashrightarrow}$ be an $\E$-triangle. For any $L\in\X$, by Lemma \ref{exact}, we have the following exact sequence
$$\C(Z,L)\xrightarrow{}\C(Y\oplus M,L)\xrightarrow{{\binom{-f}{h}}^{\ast}}\C(X,L).$$
Since $h$ is a left $\X$-approximation of $X$, for any $\beta\in \C(X,L)$, there is $n\in \C(M,L)$ such that $\beta=nh$. Take $m=(0,n)\in\C(Y\oplus M,L)$, then we have that $\beta=(0,n)\binom{-f}{h}$, that is, $X\xrightarrow{}Y\oplus M\xrightarrow{}Z\overset{}{\dashrightarrow}$ is $\C(-,\X)$-exact.

(2) It is similar to (1).\qed

\begin{lemma}\label{11l2}
Let $(\C, \E, \s)$ be an extriangulated category and $\X$ a subcategory of $\C$.
\begin{itemize}
\rm\item[(1)] If $\X$ satisfies the condition $(2^{'})$, then $\C/\X$ has cokernels.
\item[(2)] If $\X$ satisfies the condition $(3^{'})$, then $\C/\X$ has kernels.
\end{itemize}
\end{lemma}
\proof
We just prove the second statement, the first statement proves similarly. For any morphism $f\colon X\to Y$ in $\C$, there is an $\E$-triangle $$X_1\xrightarrow{}X_0\xrightarrow{\beta}Y\overset{\delta}{\dashrightarrow},$$ with $X_0,X_1\in\X$ and $\beta$ is a right $\X$-approximation of $Y$. By the dual of Lemma \ref{1.20}, we have the following commutative diagram
\begin{equation}\label{t0}
\begin{array}{l}
\xymatrix{
X_1 \ar[r]\ar@{=}[d]& K\ar[r]^{k} \ar[d]^{} & X\ar@{-->}[r]^{}\ar[d]^{f}&\\
X_1\ar[r]^{} & X_0 \ar[r]^{\beta} & Y\ar@{-->}[r]^{\del} &}
\end{array}
\end{equation}
of $\E$-triangles in $\C$. We claim that $\underline{k}$ is the kernel of $\underline{f}$. Indeed, let $g\colon U\to X$ be a morphism in $\C$ with $\underline f \underline g=0$. Then there exists $P\in\X$ such that $fg$ factor through $P$. Note that $\beta$ is a right $\X$-approximation of $Y$, $fg$ factor through $X_0$.
By \cite[Lemma 3.13]{NP}, the right square in the diagram (\ref{t0})
$$\xymatrix{K\ar[r]^{k}\ar[d]&X\ar[d]^{f}\\
X_0\ar[r]^{\beta}&Y}$$
is a weak pullback. Thus there exists a morphism $i\colon U\to K$ which makes the following
diagram
$$\xymatrix{U\ar@{-->}[dr]^{i}\ar@/_/[ddr]_w \ar@/^/[drr]^g\\&K\ar[r]^{k}\ar[d]&X\ar[d]^{f}\\
&X_0\ar[r]^{\beta}&Y}$$
commutative, we get that $\underline{k}$ is a weak kernel of $\underline{f}$.

Next, we need to show that $\underline{k}$ is a monomorphism in $\C/\X$. Suppose that $h\colon L\to K$ satisfies $\underline{k}\underline{h}=0$. Then there $U'\in\X$, $i\colon L\to U'$ and $j\colon U'\to X$ such that $kh=ji$. Since $\beta$ is a right $\X$-approximation of $Y$, there is a morphism  $l\colon U'\to X_0$ such that $fj=\beta l$. By the universal property of the weak pullback, there is a morphism $n\colon U'\to K$ such that $j=kn$. That is, we have the following commutative diagram
$$\xymatrix{ L\ar@{-->}[r]^{i}\ar[d]^h&U'\ar@{-->}[d]^{j}\ar@{-->}[dl]_n \ar@{-->}[ddl]_k \\K\ar[r]_{l}\ar[d]&X\ar[d]^{f}\\
X_0\ar[r]^{\beta}&Y}$$
in $\C$. So we have $k(ni-h)=0$. Note that the following exact sequence
$$\C(L,X_1)\xrightarrow{}\C(L,K)\xrightarrow{k_{\ast}}\C(L,X),$$
we have that the morphism $ni-h$ factors through $X_1$, hence $\underline{ni-h}=0$. It implies that $\underline{h}=0$. Therefore, $\underline{k}$ is the kernel of $\underline{f}$ in $\C/\X$.
\qed

The following lemma holds in any additive category.

\begin{lemma}{\rm\cite[Theorem 2.2]{H}\label{l2}}
Suppose that $\C$ is an additive category and $\X$ a subcategory of $\C$. Let $f\colon X\to Y$ be a morphism in $\C$. Then $\underline{f}$ is an isomorphism in $\C/\X$ if and only if there exist $P,Q\in\X$, such that the following diagram
$$\xymatrix{X\ar[r]^{f}\ar[d]^{i}&Y\\
X\oplus P\ar[r]^{\tilde{f}}&Y\oplus Q\ar[u]^{p}}$$
commutative in $\C$, where $\tilde{f}$ is an isomorphism and $i, p$ are the canonical injection and projective respectively.
\end{lemma}

\begin{lemma}\label{1l2}
Let $\C$ be an extriangulated category and $\X$ a subcategory of $\C$.
\begin{itemize}
\rm\item[(1)] If $\X$ satisfies the condition $(2^{'})$, then for any $\E$-triangle $X\xrightarrow{f}Y\xrightarrow{g}Z\overset{}{\dashrightarrow}$ which is $\C(-,\X)$-exact, we have that $\underline{g}$ is the cokernel of $\underline{f}$.
\item[(2)] If $\X$ satisfies the condition $(3^{'})$, then for any $\E$-triangle $X\xrightarrow{f}Y\xrightarrow{g}Z\overset{}{\dashrightarrow}$ which is $\C(\X,-)$-exact, we have that $\underline{f}$ is the kernel of $\underline{g}$.

\end{itemize}

\end{lemma}

\proof We just prove the first statement statement, the second statement proves similarly. For any $X\in \C$, there exists an $\E$-triangle $$X\xrightarrow{\beta}Q_0\xrightarrow{}Q_1\overset{}{\dashrightarrow},$$ with $Q_0,Q_1\in\X$ and $\beta$ is a left $\X$-approximation of $X$. By \cite[Proposition 3.15]{NP}, we have the following commutative diagram
\begin{equation}\label{tt}
\begin{array}{l}
\xymatrix{X\ar[r]^{f}\ar[d]^{\beta}&Y\ar[r]^{g}\ar[d]^{k}&Z\ar@{=}[d]\ar@{-->}[r]&\\
Q_0\ar[r]^{r}\ar[d]&C\ar[d]^{t}\ar[r]^{s}&Z\ar@{-->}[r]&\\
Q_1\ar@{-->}[d]\ar@{=}[r]&Q_1\ar@{-->}[d]&\\
&&&&}
\end{array}
\end{equation}
Let $Q\in\X$ and $i\colon Q_0\to Q$ be an morphism in $\C$. Since the $\E$-triangle $$X\xrightarrow{f}Y\xrightarrow{g}Z\overset{}{\dashrightarrow}$$ is $\C(-,\X)$-exact, we have the following exact sequence
$$\C(Z,Q)\xrightarrow{}\C(Y,Q)\xrightarrow{f^{\ast}}\C(X,Q).$$ Thus there is a morphism $j\colon Y\to Q$, such that $i\beta=jf$. Note that the upper-left square in the diagram (\ref{tt})
is a weak pushout, thus there exists a morphism $m\colon C\to Q$ which makes the following
diagram
$$\xymatrix{
X\ar[r]^f \ar[d]_{\beta} &Y \ar[d]^{k} \ar@/^/[ddr]^j\\
Q_0 \ar[r]^r \ar@/_/[drr]_i &C \ar@{-->}[dr]^m\\
&&Q}
$$
commutative. This shows that the second row is $\C(-,\X)$-exact in the diagram (\ref{tt}). Since $Q_0\in \X$, take $Q=Q_0$, we have a morphism $n\colon C\to Q$ such that $nr=1$. So we get that the second row splits in the diagram (\ref{tt}) with $C\cong Q_0\oplus Z$. Consider the following commutative diagram
$$\xymatrix{C\ar[r]^{s}\ar[d]^{id}&Z\\
C\ar[r]^{\cong\quad}&Q_0\oplus Z\ar[u]^{p}}$$
Then $\underline{s}$ is an isomorphism in $\C/\X$ by Lemma \ref{l2}. Using a dual arguments as in the proof of Lemma \ref{11l2}, we have that $\underline{k}$ is the cokernel of $\underline{f}$. Therefore, $\underline{g}$ is the cokernel of $\underline{f}$.
\qed
\vspace{3mm}

{\bf Now we give the proof of Theorem \ref{semi-abel}.}

\proof Let  $f\colon X\to Y$ be a morphism in $\C$. By Lemma \ref{l1}, we have an $\E$-triangle $$X\xrightarrow{}Y\oplus M\xrightarrow{}Z\overset{}{\dashrightarrow}$$ with $M\in\X$, which is $\C(-,\X)$-exact. By Lemma \ref{1l2}, we know that $\underline{Z}$ is the cokernel of $\underline{f}$. Applying Lemma \ref{l1} again, there is an $\E$-triangle $$X'\xrightarrow{}Y\oplus M\oplus N\xrightarrow{}Z\overset{}{\dashrightarrow}$$ with $N\in\X$, which is $\C(\X,-)$-exact. Then $\underline{X'}$ is the kernel of $c_{\underline{f}}$. By \cite[Proposition 3.17]{NP}, we obtain the following commutative diagram
\begin{equation}\label{ttt}
\begin{array}{l}
\xymatrix{X\ar[r]^{}\ar[d]^{g}&Y\oplus M\ar[r]^{\quad c_{\underline{f}}}\ar[d]^{i}&Z\ar@{=}[d]\ar@{-->}[r]&\\
X'\ar[r]^{}\ar[d]&Y\oplus M\oplus N\ar[d]^{}\ar[r]^{}&Z\ar@{-->}[r]&\\
N\ar@{-->}[d]\ar@{=}[r]&N\ar@{-->}[d]&\\
&&}
\end{array}
\end{equation}
We claim that $\underline{g}$ is an epimorphism in $\C/\X$. Indeed, for any $h\colon X\to M'$ with $M'\in \X$, since the first row is $\C(-,\X)$-exact in the diagram (\ref{ttt}), there is a morphism $j:Y\oplus M\to M'$, such that $h$ factors through $X\to Y\oplus M$. Note that the second column splits, $j$ factors through $i$. That is, we have
the following commutative diagram
\newpage
$$\xymatrix{&X\ar@{-->}[ddl]_h\ar[r]\ar[d]^{g}&Y\oplus M\ar@{-->}[ddll]^j\ar[d]^{i}\\
&X'\ar[r]^{}&Y\oplus M\oplus N\ar@{-->}[dll]\\
M'&&}
$$
This shows that $h$ factors through $g$, and hence the first column is $\C(-,\X)$-exact in the diagram (\ref{ttt}). By Lemma \ref{1l2}, $\underline{N}$ is the cokernel of $\underline{g}$. Note that $N\in\X$, we have that $\underline{g}$ is an epimorphism in $\C/\X$. Since $\C/\X$ has both kernels and cokernels, and hence we obtain the following commutative diagram
$$\xymatrix{ {\rm coker}k_{\underline{f}}\ar[r]^{\tilde{\underline{f}}}&{\rm ker}c_{\underline{f}}=\underline {X'}\ar[d]^{k_{{c_{\underline{f}}}}}  \\\underline {X}\ar[ur]^{\underline{g}}\ar[u]^{c_{{k_{\underline{f}}}}}\ar[r]^{\underline{f}}&\underline {Y}\ar[d]^{c_{\underline{f}}}\\
{\rm ker}{\underline{f}}\ar[u]^{k_{\underline{f}}}&{\rm coker}{\underline{f}}=\underline{Z}}$$
in $\C/\X$. So the canonical morphism $\tilde{\underline{f}}$ is an epimorphism. Dually, we can prove $\tilde{\underline{f}}$ is a monomorphism. Therefore, the quotient category $\C/\X$ is a semi-abelian category.
\qed
\vspace{2mm}

By applying Theorem \ref{semi-abel} to exact categories, and using the fact that any exact
category can be viewed as an extriangulated category, we get the following result.

\begin{corollary}\rm\cite[Theorem 3.4]{XZ}
Let $\C$ be exact category, and $\X$ a pseudo cluster tilting subcategory of $\C$. Then the quotient category $\C/\X$ is a semi-abelian category.
\end{corollary}
\proof It follows from Remark \ref{exam} and Remark \ref{epst}.
\qed
\vspace{2mm}

In order to consider when the quotient category $\C/\X$ is abelian category, we introduce two special classes of $\E$-triangles $\S$ and $\T$.

$\bullet$ An $\E$-triangle $S_1\overset{}{\longrightarrow}S_2\overset{}{\longrightarrow}S_3\overset{}{\dashrightarrow}$ belongs to $\S$, if there exists a commutative diagram of  $\E$-triangles
$$\xymatrix{U\ar[r]^{}\ar@{=}[d]&S_1\ar[r]^{}\ar[d]^{}&M_1\ar[d]\ar@{-->}[r]&\\
U\ar[r]^{}&S_2\ar[d]^{}\ar[r]^{}&M_2\ar[d]^{}\ar@{-->}[r]&\\
&S_3\ar@{-->}[d]\ar@{=}[r]&S_3\ar@{-->}[d]\\
&&&}$$
where the $\E$-triangle $$U\overset{}{\longrightarrow}S_2\overset{}{\longrightarrow}M_2\overset{}{\dashrightarrow}$$ is $\C(\X,-)$-exact, and the $\E$-triangle $$M_1\overset{}{\longrightarrow}M_2\overset{}{\longrightarrow}S_3\overset{}{\dashrightarrow}$$ is $\C(-,\X)$-exact.

$\bullet$ Dually, an $\E$-triangle $T_1\overset{}{\longrightarrow}T_2\overset{}{\longrightarrow}T_3\overset{}{\dashrightarrow}$ belongs to $\T$, if there exists a commutative diagram of  $\E$-triangles
$$\xymatrix{T_1\ar[r]^{}\ar@{=}[d]&N_2\ar[r]^{}\ar[d]^{}&N_3\ar[d]^{}\ar@{-->}[r]&\\
T_1\ar[r]^{}&T_2\ar[d]^{}\ar[r]^{}&T_3\ar[d]^{}\ar@{-->}[r]&\\
&V\ar@{-->}[d]\ar@{=}[r]&V\ar@{-->}[d]\\
&&&}$$
where the $\E$-triangle $$N_2\overset{}{\longrightarrow}T_2\overset{}{\longrightarrow}V\overset{}{\dashrightarrow}$$ is $\C(-,\X)$-exact, and the $\E$-triangle $$T_1\overset{}{\longrightarrow}N_2\overset{}{\longrightarrow}N_3\overset{}{\dashrightarrow}$$ is $\C(\X,-)$-exact.

Let $\Y$ be a class of $\E$-triangles in $\C$. If any $\E$-triangle $D\overset{}{\longrightarrow}E\overset{}{\longrightarrow}F\overset{}{\dashrightarrow}$ with $D, F\in \X$ belonging to $\Y$ splits, then we say that $\X$ is  \emph{self-orthogonal} with respect to $\Y$.

Our second main result is the following.
\begin{theorem}\label{abel}
Let $\C$ be extriangulated category and $\X$ a pseudo cluster tilting subcategory of $\C$. Then the following statements are
equivalent.
\begin{itemize}
 \item[\rm (1)] $\C/\X$ is an abelian category.
\item[\rm (2)] $\X$ is self-orthogonal with respect to $\S$.
\item[\rm (3)] $\X$ is self-orthogonal with respect to $\T$.

\end{itemize}
\end{theorem}
\proof
$(1)\Rightarrow(2)$ Suppose that $S_1\overset{}{\longrightarrow}S_2\overset{}{\longrightarrow}S_3\overset{}{\dashrightarrow}$ is an $\E$-triangle in $\S$, where $S_1,S_3\in \X$. Then there exists a commutative diagram of $\E$-triangles
$$\xymatrix{U\ar[r]^{}\ar@{=}[d]&S_1\ar[r]^{}\ar[d]^{}&M_1\ar[d]^{m}\ar@{-->}[r]&\\
U\ar[r]^{}&S_2\ar[d]^{p}\ar[r]^{r}&M_2\ar[d]^{q}\ar@{-->}[r]&\\
&S_3\ar@{-->}[d]\ar@{=}[r]&S_3\ar@{-->}[d]\\
&&&}$$
where the $\E$-triangle $M_1\overset{m}{\longrightarrow}M_2\overset{q}{\longrightarrow}S_3\overset{}{\dashrightarrow}$ is $\C(-,\X)$-exact. By Lemma \ref{1l2}, we have that $\underline{m}$ is an epimorphism. Since the $\E$-triangle $U\overset{}{\longrightarrow}S_2\overset{r}{\longrightarrow}M_2\overset{}{\dashrightarrow}$ is $\C(\X,-)$-exact, one
can check that $\underline{S_1}$ is a weak kernel of $\underline{m}$, and hence $\underline{m}$ is a monomorphism. Hence $\underline{m}$ is an isomorphism since $\C/\X$ is an abelian category. We obtain the following commutative diagram
$$\xymatrix{M_1\ar[r]^{m}\ar[d]^{i}&M_2\\
M_1\oplus P\ar[r]^{\tilde{m}}&M_2\oplus Q\ar[u]^{p}}$$
where $\tilde{m}$ is an isomorphism, $P,Q\in \X$. Set $\tilde{m} =\left(
             \begin{smallmatrix}
              m & b \\
               c & d \\
             \end{smallmatrix}
           \right)$ and  $\tilde{m}^{-1} =\left(
             \begin{smallmatrix}
              m_1 & b_1 \\
               c_1 & d_1 \\
             \end{smallmatrix}
           \right)$.
Since the $\E$-triangle $M_1\overset{m}{\longrightarrow}M_2\overset{q}{\longrightarrow}S_3\overset{}{\dashrightarrow}$ is $\C(-,\X)$-exact, there exists $v:M_2\to Q$ such that $c=vm$. Note that $m_1m+b_1c=1$, we have $(m_1+b_1v)m=1$. Hence $m$ is a section, and the $\E$-triangle $M_1\overset{m}{\longrightarrow}M_2\overset{q}{\longrightarrow}S_3\overset{}{\dashrightarrow}$ is split. This means the $\E$-triangle $S_1\overset{}{\longrightarrow}S_2\overset{}{\longrightarrow}S_3\overset{}{\dashrightarrow}$ is also split. So $\X$ is self-orthogonal with respect to $\S$.

$(2)\Rightarrow(1)$ We need to show that a regular morphism in $\C/\X$ must be an isomorphism. We use the same notation as the proof of Theorem \ref{semi-abel}. Let $\underline {f}\colon \underline {X}\to \underline {Y}$ be a monomorphism in $\C/\X$ induced by $f\colon X\to Y$ in $\C$. We obtain an $\E$-triangle $X\overset{g}{\longrightarrow}X'\overset{}{\longrightarrow}N\overset{}{\dashrightarrow}$ with $N\in \X$, which is $\C(-,\X)$-exact as in the diagram (\ref{ttt}). For the object $X'$, there exists an $\E$-triangle $P_1\overset{}{\longrightarrow}P_0\overset{}{\longrightarrow}X'\overset{}{\dashrightarrow}$ with $P_0,P_1\in \X$, which is $\C(\X,-)$-exact. By (ET4)$\op$, we have the following commutative diagram of $\E$-triangles
$$\xymatrix{P_1\ar[r]^{}\ar@{=}[d]&K\ar[r]^{}\ar[d]^{}&X\ar[d]^{g}\ar@{-->}[r]&\\
P_1\ar[r]^{}&P_0\ar[d]^{}\ar[r]^{}&X'\ar[d]^{}\ar@{-->}[r]&\\
&N\ar@{-->}[d]\ar@{=}[r]&N\ar@{-->}[d]\\
&&&}$$
This show that the $\E$-triangle $K\overset{}{\longrightarrow}P_0\overset{}{\longrightarrow}N\overset{}{\dashrightarrow}$ belongs to $\S$. Consider the following commutative diagram
$$\xymatrix{X\ar[r]^{g}\ar[d]^{}&X'\ar[d]^{}\\
Y\oplus M\ar[r]^{i}&Y\oplus M\oplus N}$$
in $\C$ from the diagram (\ref{ttt}), which induces the following  commutative diagram
\begin{equation}\label{mmm}
\begin{array}{l}\xymatrix{\underline {X}\ar[r]^{\underline {g}}\ar[d]^{\underline {f}}&\underline {X'}\ar[d]^{k_{c_{\underline {f}}}}\\
\underline {Y}\oplus \underline {M}\ar[r]^{\underline{i}}&\underline {Y}\oplus \underline {M}\oplus \underline {N}}
\end{array}
\end{equation}
in $\C/\X$. Since $\underline{i}$ is an isomorphism by Lemma \ref{l2} and $\underline {f}$ is a monomorphism, we know that $\underline {g}$ is a monomorphism. So $K\in\X$, since $\underline{K}$ is the kernel of $\underline {g}$ by Lemma \ref{11l2}. We have that the $\E$-triangle $K\overset{}{\longrightarrow}P_0\overset{}{\longrightarrow}N\overset{}{\dashrightarrow}$ splits. This means the $\E$-triangle $X\overset{g}{\longrightarrow}X'\overset{}{\longrightarrow}N\overset{}{\dashrightarrow}$ is also split. Hence in the diagram \ref{mmm}, the two rows are all isomorphisms in $\C/\X$. We get that $\underline{f}$ is a kernel. That means that any monomorphism in $\C/\X$ is a kernel. The rest is well-known and we leave it to the reader.

$(1)\Longleftrightarrow(3)$ The equivalence is the dual of $(1)\Longleftrightarrow(2)$.  \qed

\vspace{2mm}

By applying Theorem \ref{semi-abel} to exact categories, and using the fact that any exact
category can be viewed as an extriangulated category, we get the following result.

\begin{corollary}
In Theorem {\rm \ref{abel}}, when $\C$ is an exact category, it is just Theorem {\rm 3.4} in {\rm \cite{XZ}}.
\end{corollary}

We provide some examples to explain our main results.

\begin{example}
(1) Suppose that the category $\C$ is as shown in Example \ref{00ex}, we know that the subcategory $\X:=\add(2\oplus\begin{aligned}1\\[-3mm]2\end{aligned}\oplus \begin{aligned}2\\[-3mm]1\\[-3mm]2\end{aligned})$ is a pseudo cluster tilting subcategory of $\C$. So the quotient category $\C/\X$ is a semi-abelian category by Theorem \ref{semi-abel}.
\vspace{2mm}

(2) Suppose that the category $\C$ is as shown in Example \ref{ex2}, we know that the subcategory
$\Y=\add(2\oplus1\oplus\begin{aligned}2\\[-3mm]1\end{aligned}\oplus 2[2]\oplus\begin{aligned}1\\[-3mm]2\end{aligned}[2])$ is a pseudo cluster tilting subcategory of $\C$. So the quotient category $\C/\Y$ is a semi-abelian category by Theorem \ref{semi-abel}.

\end{example}

\textbf{Jian He}\\
Department of Applied Mathematics, Lanzhou University of Technology, 730050 Lanzhou, Gansu, P. R. China\\
E-mail: \textsf{jianhe30@163.com}\\[0.3cm]
\textbf{Jing He}\\
College of Science, Hunan University of Technology and Business, 410205 Changsha, Hunan P. R. China\\
E-mail: \textsf{jinghe1003@163.com}

\end{document}